\documentclass[psamsfonts]{amsart}

\usepackage{amssymb,amsfonts}
\usepackage[all,arc]{xy}
\usepackage{enumerate}
\usepackage{mathrsfs}
\usepackage{graphicx}

\newtheorem{theorem}{Theorem}[section]
\newtheorem{corollary}[theorem]{Corollary}
\newtheorem{proposition}[theorem]{Proposition}
\newtheorem{lemma}[theorem]{Lemma}

\theoremstyle{definition}
\newtheorem{definition}[theorem]{Definition}

\newtheorem{question}[theorem]{Question}
\newtheorem{remark}[theorem]{Remark}

\DeclareMathOperator{\Ann}{Ann}

\DeclareMathOperator{\Id}{Id}

\DeclareMathOperator{\Hom}{Hom}


\numberwithin{equation}{section}

\begin{document}

\title[Auslander Modules]{Auslander Modules}

\author{Peyman Nasehpour}

\dedicatory{Dedicated to my father, Maestro Nasrollah Nasehpour}

\keywords{Auslander modules, Auslander's zero-divisor conjecture, content algebras, torsion-free modules}

\subjclass[2010]{13A15, 13B25, 13F25}

\address{Department of Engineering Science, Golpayegan University of Technology, Golpayegan, Iran}
\email{nasehpour@gut.ac.ir, nasehpour@gmail.com}

\maketitle

\begin{abstract}

In this paper, we introduce the notion of Auslander modules, inspired from Auslander's zero-divisor conjecture (theorem) and give some interesting results for these modules. We also investigate torsion-free modules.

\end{abstract}

\section{Introduction}

Auslander's zero-divisor conjecture in commutative algebra states that if $R$ is a Noetherian local ring, $M$ is a nonzero $R$-module of finite type, and finite projective dimension, and $r\in R$ is not a zero-divisor on $M$, then $r$ is not a zero-divisor on $R$ \cite[p. 8]{Hochster1975} and \cite[p. 496]{Hochster1987}. This ``conjecture'' is in fact a theorem, after Peskine and Szpiro in \cite{PeskineSzpiro1973} showed that Auslander's zero-divisor theorem is a corollary of their new intersection theorem and thereby proved it for a large class of local rings. Also see \cite[p. 417]{Roberts1989}. Note that its validity without any restrictions followed when Roberts \cite{Roberts1987} proved the new intersection theorem  in full generality. Also see Remark 9.4.8 in \cite{BrunsHerzog1998}.

Let $M$ be an arbitrary unital nonzero module over a commutative ring $R$ with a nonzero identity. Inspired from Auslander's zero-divisor theorem, one may ask when the inclusion $Z_R(R) \subseteq Z_R(M)$ holds, where by $Z_R(M)$, we mean the set of all zero-divisors of the $R$-module $M$. In Definition \ref{Auslandermodule}, we define an $R$-module $M$ to be Auslander if $Z_R(R) \subseteq Z_R(M)$ and in Proposition \ref{AuslandermoduleEx}, we give a couple of examples for the families of Auslander modules. The main theme of $\S$1 is to see under what conditions if $M$ is an Auslander $R$-module, then the $S$-module $M \otimes_R S$ is Auslander, where $S$ is an $R$-algebra (see Theorem \ref{AuslandermoduleThm}, Theorem \ref{FormalPowerAuslander}, and Theorem \ref{AuslanderMcCoy}). For example, in Corollary \ref{AuslanderContent}, we show that if $M$ is an Auslander $R$-module, $B$ a content $R$-algebra, and $M$ has property (A), then $M \otimes_R B$ is an Auslander $B$-module. For
the definition of content algebras refer to \cite[Section 6]{OhmRush1972}.

On the other hand, let us recall that an $R$-module $M$ is torsion-free if the natural map $M \rightarrow M \otimes Q$ is injective, where $Q$ is the total quotient ring of the ring $R$ \cite[p. 19]{BrunsHerzog1998}. It is easy to see that $M$ is a torsion-free $R$-module if and only if $Z_R(M) \subseteq Z_R(R)$. In $\S$2, we investigate torsion-free property under polynomial and power series extensions (see Theorem \ref{PropertyZThm} and Theorem \ref{FormalPowerZ}). We also investigate torsion-free Auslander modules (check Proposition \ref{Z-AuslanderPro}, Theorem \ref{torsionfreeAuslanderMcCoy}, and Theorem \ref{Z-AuslanderThm}).

In this paper, all rings are commutative with non-zero identities and all modules are unital.

\section{Auslander Modules}

We start the first section by defining Auslander modules:

\begin{definition}

\label{Auslandermodule}

We define an $R$-module $M$ to be an Auslander module, if $r\in R$ is not a zero-divisor on $M$, then $r$ is not a zero-divisor on $R$, or equivalently, if the following property holds: $$Z_R(R) \subseteq Z_R(M).$$
\end{definition}

Let us recall that if $M$ is an $R$-module, the content of $m\in M$, denoted by $c(m)$, is defined to be the following ideal: $$c(m) = \bigcap \{I\in \Id(R): m\in IM\},$$ where by $\Id(R)$, we mean the set of all ideals of $R$. The $R$-module $M$ is said to be a content $R$-module, if $m\in c(m)M$, for all $m\in M$ \cite{OhmRush1972}. In the following, we give some families of Auslander modules:

\begin{proposition}[Some Families of Auslander Modules] Let $M$ be an $R$-module. Then the following statements hold:
\label{AuslandermoduleEx}
\begin{enumerate}

\item If $R$ is a domain, then $M$ is an Auslander $R$-module.

\item If $M$ is a flat and content $R$-module such that for any $s\in R$, there is an $x\in M$ such that $c(x) = (s)$. Then $M$ is an Auslander $R$-module.

\item If $M$ is an $R$-module such that $\Ann(M)=(0)$, then $M$ is an Auslander $R$-module.

\item If for any nonzero $s\in R$, there is an $x\in M$ such that $s \cdot x \neq 0$, i.e. $\Ann(M) =(0)$, then $\Hom_R(M,M)$ is an Auslander $R$-module.

\item If $N$ is an $R$-submodule of an $R$-module $M$ and $N$ is Auslander, then $M$ is also Auslander.

\item If $M$ is an Auslander $R$-module, then $M \oplus M^{\prime}$ is an Auslander $R$-module for any $R$-module $M^{\prime}$. In particular, if $\{M_i\}_{i\in \Lambda}$ is a family of $R$-modules and there is an $i\in \Lambda$, say $i_0$, such that $M_{i_0}$ is an Auslander $R$-module, then $\bigoplus_{i\in \Lambda} M_i$ and $\prod_{i\in \Lambda} M_i$ are Auslander $R$-modules.

\end{enumerate}

\begin{proof}

The statement $(1)$ is obvious. We prove the other statements:

$(2)$: Let $r\in Z_R(R)$. By definition, there is a nonzero $s\in R$ such that $r\cdot s = 0$. Since in content modules $c(x) = (0)$ if and only if $x = 0$ \cite[Statement 1.2]{OhmRush1972} and by assumption, there is a nonzero $x\in M$ such that $c(x) = (s)$, we obtain that $r \cdot c(x) = (0)$. Also, since $M$ is a flat and content $R$-module, by \cite[Theorem 1.5]{OhmRush1972}, $r\cdot c(x) = c(r \cdot x)$. This implies that $r\in Z_R(M)$.

$(3)$: Suppose that $r \cdot s = 0$ for some nonzero $s$ in $R$. By assumpstion, there exists an $x$ in $M$ such that $s \cdot x \neq 0$, but $r\cdot (s\cdot x) = 0$, and so $r$ is a zero-divisor on $M$.

$(4)$: Let $r\in Z_R(R)$. So, there is a nonzero $s\in R$ such that $r \cdot s=0$. Define $f_s : M \longrightarrow M$ by $f_s(x) = s \cdot x$. By assumption, $f_s$ is a nonzero element of $\Hom_R(M,M)$. But $rf_s = 0$. This means that $r\in Z_R(\Hom(M,M))$.

$(5)$: The proof is straightforward, if we consider that $Z(N) \subseteq Z(M)$.

The statement $(6)$ is just a corollary of the statement $(5)$.
\end{proof}

\end{proposition}

\begin{proposition}
Let $M$ be an Auslander $R$-module and $S$ a multiplicatively closed subset of $R$ contained in $R-Z_R(M)$. Then, $M_S$ is an Auslander $R_S$-module.

\begin{proof}
Let $Z_R(R) \subseteq Z_R(M)$ and $S$ be a multiplicatively closed subset of $R$ such that $S \subseteq R-Z_R(M)$. Take $r_1 / s_1 \in Z_{R_S} (R_S)$. So there exists an $r_2 / s_2 \neq 0/1$ such that $(r_1 \cdot r_2) / (s_1 \cdot s_2) = 0/1$. Since $S \subseteq R-Z_R(R)$, we have $r_1 \cdot r_2 = 0$, where $r_2 \neq 0$. But $Z_R(R) \subseteq Z_R(M)$, so $r_1 \in Z_R(M)$. Consequently, there is a nonzero $m\in M$ such that $r_1 \cdot m = 0$. Since $S \subseteq R-Z_R(M)$, $m/1$ is a nonzero element of $M_S$. This point that $r_1 / s_1 \cdot m / 1 = 0 / 1$, causes $r_1 / s_1$ to be an element of $Z_{R_S} (M_S)$ and the proof is complete.
\end{proof}

\end{proposition}

Let us recall that an $R$-module $M$ has property (A), if each finitely generated ideal $I \subseteq Z_R(M)$ has a nonzero annihilator in $M$ \cite[Definition 10]{Nasehpour2011}. Examples of modules having property (A) include modules having very few zero-divisors \cite[Definition 6]{Nasehpour2011}. Especially, finitely generated modules over Noetherian rings have property (A) \cite[p. 55]{Kaplansky1970}. Homological aspects of modules having very few zero-divisors have been investigated in \cite{NashpourPayrovi2010}. Finally, we recall that if $R$ is a ring, $G$ a monoid, and $f = r_1 X^{g_1} + \cdots + r_n X^{g_n}$ is an element of the monoid ring $R[G]$, then the content of $f$, denoted by $c(f)$, is the finitely generated ideal $(r_1, \ldots, r_n)$ of $R$.

\begin{theorem}

\label{AuslandermoduleThm}

Let the $R$-module $M$ have property (A) and $G$ be a commutative, cancellative, and torsion-free monoid. Then, $M[G]$ is an Auslander $R[G]$-module if and only if $M$ is an Auslander $R$-module.
\begin{proof}

$(\Rightarrow)$: Let $r\in Z_R(R)$. So, $r\in Z_{R[G]}(R[G])$ and by assumption, $r\in Z_{R[G]}(M[G])$. Clearly, this means that there is a nonzero $g$ in $Z_{R[G]}(M[G])$ such that $rg =0$. Therefore, there is a nonzero $m$ in $M$ such that $rm =0$.

$(\Leftarrow)$: Let $f\in Z_{R[G]}(R[G])$. By \cite[Theorem 2]{Nasehpour2011}, there is a nonzero element $r\in R$ such that $f\cdot r = 0$. This implies that $c(f) \subseteq Z_R(R)$. But $M$ is an Auslander $R$-module, so $Z_R(R) \subseteq Z_R(M)$, which implies that $c(f) \subseteq Z_R(M)$. On the other hand, $M$ has property (A). Therefore, $c(f)$ has a nonzero annihilator in $M$. Hence, $f\in Z_{R[G]}(M[G])$ and the proof is complete.
\end{proof}

\end{theorem}

Note that a semimodule version of Theorem \ref{AuslandermoduleThm} has been given in \cite{Nasehpour2017}.

It is good to mention that if $R$ is a ring and $f = a_0 + a_1 X + \cdots + a_n X^n + \cdots$ is an element of $R[[X]]$, then $A_f$ is defined to be the ideal of $R$ generated by the coefficients of $f$, i.e. $$A_f = (a_0 , a_1 , \ldots , a_n , \ldots).$$ One can easily check that if $R$ is Noetherian, then $A_f = c(f)$. The following lemma is a generalization of Theorem 5 in \cite{Fields1971}:

\begin{lemma}

\label{McCoyPowerSeries}

Let $R$ be a Noetherian ring, $M$ a finitely generated $R$-module, $f\in R[[X]]$, $g\in M[[X]]-\{0\}$, and $fg = 0$. Then, there is a nonzero constant $m\in M$ such that $f \cdot m = 0$.

\begin{proof}
Define $c(g)$, the content of $g$, to be the $R$-submodule of $M$ generated by its coefficients. If $c(f)c(g)=(0)$, then choose a nonzero $m\in c(g)$. Clearly, $f\cdot m = 0$. Otherwise, by Theorem 3.1 in \cite{EpsteinShapiro2016AMS}, one can choose a positive integer $k$, such that $c(f) c(f)^{k-1} c(g) = 0$, while $c(f)^{k-1} c(g) \neq 0$. Now for each nonzero element $m$ in $c(f)^{k-1} c(g)$, we have $f\cdot m = 0$ and the proof is complete.
\end{proof}
\end{lemma}

\begin{theorem}

\label{FormalPowerAuslander}

Let $R$ be a Noetherian ring and the $R$-module $M$ have property (A). Then, $M[[X]]$ is an Auslander $R[[X]]$-module if and only if $M$ is an Auslander $R$-module.

\begin{proof}

By Lemma \ref{McCoyPowerSeries}, the proof is just a mimicking of the proof of Theorem \ref{AuslandermoduleThm}.
\end{proof}
\end{theorem}

Since finitely generated modules over Noetherian rings have property (A) \cite[p. 55]{Kaplansky1970}, we have the following corollary:

\begin{corollary}
Let $R$ be a Noetherian ring and $M$ be a finitely generated $R$-module. Then, $M[[X]]$ is an Auslander $R[[X]]$-module if and only if $M$ is an Auslander $R$-module.
\end{corollary}

\begin{remark}[Ohm-Rush Algebras] \label{ContentFunction} Let us recall that if $B$ is an $R$-algebra, then $B$ is said to be an Ohm-Rush $R$-algebra, if $f\in c(f)B$, for all $f\in B$ \cite[Definition 2.1]{EpsteinShapiro2016}. It is easy to see that if $P$ is a projective $R$-algebra, then $P$ is an Ohm-Rush $R$-algebra \cite[Corollary 1.4]{OhmRush1972}. Note that if $R$ is a Noetherian ring and $f = a_0 + a_1 X + \cdots + a_n X^n + \cdots$ is an element of $R[[X]]$, then $A_f = c(f)$, where $A_f$ is the ideal of $R$ generated by the coefficients of $f$. This simply implies that $R[[X]]$ is an Ohm-Rush $R$-algebra.

\end{remark}

Now we go further to define McCoy algebras, though we don't go through them deeply in this paper. McCoy semialgebras (and algebras) and their properties have been discussed in more details in author's recent paper on zero-divisors of semimodules and semialgebras \cite{Nasehpour2017}.

\begin{definition}
We say that $B$ is a McCoy $R$-algebra, if $B$ is an Ohm-Rush $R$-algebra and $f \cdot g = 0$ with $g \neq 0$ implies that there is a nonzero $r\in R$ such that $c(f)\cdot r = (0)$, for all $f,g \in B$.
\end{definition}

Since any content algebra is a McCoy algebra \cite[Statement 6.1]{OhmRush1972}, we have plenty of examples for McCoy algebras. For instance, if $G$ is a torsion-free abelian group and $R$ is a ring, then $R[G]$ is a content - and therefore, a McCoy - $R$-algebra \cite{Northcott1959}.  For other examples of McCoy algebras, one can refer to content algebras given in Examples 6.3 in \cite{OhmRush1972}. Now we proceed to give the following general theorem on Auslander modules:

\begin{theorem}

\label{AuslanderMcCoy}

Let $M$ be an Auslander $R$-module and $B$ a faithfully flat McCoy $R$-algebra. If $M$ has property (A), then $M \otimes_R B$ is an Auslander $B$-module.

\begin{proof}

Let $f\in Z_B(B)$. So by definition, there is a nonzero $r\in R$ such that $c(f) \cdot r = (0)$. This implies that $c(f) \subseteq Z_R(R)$. But $M$ is an Auslander $R$-module. Therefore, $c(f) \subseteq Z_R(M)$. Since $c(f)$ is a finitely generated ideal of $R$ \cite[p. 3]{OhmRush1972} and $M$ has property (A), there is a nonzero $m\in M$ such that $c(f)\cdot m = (0)$. This means that $c(f) \subseteq \Ann_R(m)$. Therefore, $c(f)B \subseteq \Ann_R(m)B$. Since any McCoy $R$-algebra is by definition an Ohm-Rush $R$-algebra, we have that $f \in c(f)B$. Our claim is that $\Ann_R(m)B = \Ann_B(1 \otimes m)$ and here is the proof: Since $$ 0 \longrightarrow R/\Ann_R(m) \longrightarrow M$$ is an $R$-exact sequence and $B$ is a faithfully flat $R$-module, we have the following $B$-exact sequence: $$ 0 \longrightarrow B/\Ann_R(m)B \longrightarrow M \otimes_R B,$$ with $\Ann_R(m)B = \Ann (m \otimes_R 1_B)$. This means that $f\in Z_B(M \otimes_R B)$ and the proof is complete.
\end{proof}

\end{theorem}

\begin{corollary}

\label{AuslanderContent}

Let $M$ be an Auslander $R$-module and $B$ a content $R$-algebra. If $M$ has property (A), then $M \otimes_R B$ is an Auslander $B$-module.
\begin{proof}
By definition of content algebras \cite[Section 6]{OhmRush1972}, any content $R$-algebra is faithfully flat. Also, by \cite[Statement 6.1]{OhmRush1972}, any content $R$-algebra is a McCoy $R$-algebra.
\end{proof}

\end{corollary}

\begin{question}
Is there any faithfully flat McCoy algebra that is not a content algebra?
\end{question}

\section{Torsion-Free Modules}

Let us recall that if $R$ is a ring, $M$ an $R$-module, and $Q$ the total ring of fractions of $R$, then $M$ is torsion-free if the natural map $M \rightarrow M \otimes Q$ is injective \cite[p. 19]{BrunsHerzog1998}. It is starightforward to see that $M$ is a torsion-free $R$-module if and only if $Z_R(M) \subseteq Z_R(R).$ Therefore, the notion of Auslander modules defined in Definition \ref{Auslandermodule} is a kind of dual to the notion of torsion-free modules.

The proof of the following theorem is quite similar to the proof of Proposition \ref{AuslandermoduleThm}. Therefore, we just mention the proof briefly.

\begin{theorem}

\label{PropertyZThm}

Let the ring $R$ have property (A) and $G$ be a commutative, cancellative, and torsion-free monoid. Then, the $R[G]$-module $M[G]$ is torsion-free if and only if the $R$-module $M$ is torsion-free.

\begin{proof}

$(\Rightarrow)$: Let $r\in Z_R(M)$. Clearly, this implies that $r\in Z_{R[G]}(M[G])$. But the $R[G]$-module $M[G]$ is torsion-free. Therefore, $Z_{R[G]}(M[G]) \subseteq Z_{R[G]}(R[G])$. So, $r\in Z_R(R)$.

$(\Leftarrow)$: Let $f\in Z_{R[G]}(M[G])$. By \cite[Theorem 2]{Nasehpour2011}, there is a nonzero $m\in M$ such that $c(f) \cdot m=0$, which means that $c(f) \subseteq Z_R(M)$. Since $M$ is torsion-free, $c(f) \subseteq Z_R(R)$, and since $R$ has property (A), $f\in Z_{R[G]}(R[G])$ and the proof is complete.
\end{proof}

\end{theorem}

\begin{theorem}

\label{FormalPowerZ}
Let $R$ be a Noetherian ring and $M$ be a finitely generated $R$-module. Then, the $R[[X]]$-module $M[[X]]$ is torsion-free if and only if the $R$-module $M$ is torsion-free.

\begin{proof}

$(\Rightarrow)$: Its proof is similar to the proof of Theorem \ref{PropertyZThm} and therefore, we don't bring it here.

$(\Leftarrow)$: Let $f\in Z_{R[[X]]}(M[[X]])$. By Lemma \ref{McCoyPowerSeries}, there is a nonzero element $m\in M$ such that $f\cdot m = 0$. By Remark \ref{ContentFunction}, this implies that $c(f) \subseteq Z_R(M)$. But $M$ is torsion-free, so $Z_R(M) \subseteq Z_R(R)$, which implies that $c(f) \subseteq Z_R(R)$. On the other hand, since every Noetherian ring has property (A) (check \cite[Theorem 82, p. 56]{Kaplansky1970}), $c(f)$ has a nonzero annihilator in $R$. This means that $f\in Z_{R[[X]]}(R[[X]])$, Q.E.D.
\end{proof}
\end{theorem}

We continue this section by investigating torsion-free Auslander modules.

\begin{remark}

In the following, we show that there are examples of modules that are Auslander but not torsion-free and also there are some modules that are torsion-free but not Auslander.

\begin{enumerate}

\item Let $R$ be a ring and $S \subseteq R-Z_R(R)$ a multiplicatively closed subset of $R$. Then, it is easy to see that $Z_R(R) = Z_R(R_S)$, i.e. $R_S$ is a torsion-free Auslander $R$-module.

\item Let $D$ be a domain and $M$ a $D$-module such that $Z_D(M) \neq \{0\}$. Clearly, $Z_D(D)=\{0\}$ and therefore, $M$ is Auslander, while $M$ is not torsion-free. For example, if $D$ is a domain that is not a field, then $D$ has an ideal $I$ such that $I\neq (0)$ and $I \neq D$. It is clear that $Z_D(D/I)\supseteq I$.

\item Let $k$ be a field and consider the ideal $I=(0) \oplus k$ of the ring $R= k \oplus k$. It is easy to see that $Z_R(R) = ((0)\oplus k) \cup (k \oplus (0))$, while $Z_R(R/I) = (0) \oplus k$. This means that the $R$-module $R/I$ is torsion-free, while it is not Auslander.

\end{enumerate}

\end{remark}

\begin{proposition}[Some Families of Torsion-free Auslander Modules] Let $M$ be an $R$-module. Then, the following statements hold:

\begin{enumerate}

\item If $R$ is a domain and $M$ is a flat $R$-module, then $M$ is torsion-free Auslander $R$-module.

\item If $M$ is a flat and content $R$-module such that for any $s\in R$, there is an $x\in M$ such that $c(x) = (s)$. Then $M$ is a torsion-free Auslander $R$-module.

\item If $R$ is a Noetherian ring and $M$ is a finitely generated flat $R$-module and for any nonzero $s\in R$, there is an $x\in M$ such that $s \cdot x \neq 0$. Then $\Hom_R(M,M)$ is a torsion-free Auslander $R$-module.

\item If $M$ is an Auslander $R$-module, and $M$ and $M^{\prime}$ are both flat modules, then $M \oplus M^{\prime}$ is a torsion-free Auslander $R$-module. In particular, if $\{M_i\}_{i\in \Lambda}$ is a family of flat $R$-modules and there is an $i\in \Lambda$, say $i_0$, such that $M_{i_0}$ is an Auslander $R$-module, then $\bigoplus_{i\in \Lambda} M_i$ is a torsion-free Auslander $R$-module.

\item If $R$ is a coherent ring and $\{M_i\}_{i\in \Lambda}$ is a family of flat $R$-modules and there is an $i\in \Lambda$, say $i_0$, such that $M_{i_0}$ is an Auslander $R$-module, then $\prod_{i\in \Lambda} M_i$ is a torsion-free Auslander $R$-module.

\end{enumerate}

\begin{proof}
It is trivial that every flat module is torsion-free. By considering Proposition \ref{AuslandermoduleEx}, the proof of statements $(1)$ and $(2)$ is straightforward.

The proof of statement $(3)$ is based on Theorem 7.10 in \cite{Matsumura1989} that says that each finitely generated flat module over a local ring is free. Now if $R$ is a Noetherian ring and $M$ is a flat and finitely generated $R$-module, then $M$ is a locally free $R$-module. This causes $\Hom_R(M,M)$ to be also a locally free $R$-module and therefore, $\Hom_R(M,M)$ is $R$-flat and by Proposition \ref{AuslandermoduleEx}, a torsion-free Auslander $R$-module.

The proof of the statements $(4)$ and $(5)$ is also easy, if we note that the direct sum of flat modules is flat \cite[Proposition 4.2]{Lam1999} and, if $R$ is a coherent ring, then the direct product of flat modules is flat \cite[Theorem 4.47]{Lam1999}.
\end{proof}

\end{proposition}

\begin{proposition}

\label{Z-AuslanderPro}

Let both the ring $R$ and the $R$-module $M$ have property (A) and $G$ be a commutative, cancellative, and torsion-free monoid. Then, $M[G]$ is a torsion-free Auslander $R[G]$-module if and only if $M$ is a torsion-free Auslander $R$-module.

\begin{proof}
By Theorem \ref{AuslandermoduleThm} and Theorem \ref{PropertyZThm}, the statement holds.
\end{proof}
\end{proposition}

\begin{corollary}
Let $R$ be a Noetherian ring and $M$ a finitely generated $R$-module, and $G$ a commutative, cancellative, and torsion-free monoid. Then, $M[G]$ is a torsion-free Auslander $R[G]$-module if and only if $M$ is a torsion-free Auslander $R$-module.
\end{corollary}

\begin{theorem}
\label{torsionfreeAuslanderMcCoy}
Let $M$ be a flat Auslander $R$-module and $B$ a faithfully flat McCoy $R$-algebra. If $M$ has property (A), then $M \otimes_R B$ is a torsion-free Auslander $B$-module.

\begin{proof}
By Theorem \ref{AuslanderMcCoy}, $Z_B(B) \subseteq Z_B(M \otimes_R B)$. On the other hand, since $M$ is a flat $R$-module, by \cite[Proposition 4.1]{Lam1999}, $M \otimes_R B$ is a flat $B$-module. This implies that $Z_B(M \otimes_R B) \subseteq Z_B(B)$ and the proof is complete.
\end{proof}
\end{theorem}

\begin{corollary}
Let $M$ be a flat Auslander $R$-module and $B$ a content $R$-algebra. If $M$ has property (A), then $M \otimes_R B$ is a torsion-free Auslander $B$-module.
\end{corollary}

\begin{theorem}

\label{Z-AuslanderThm}

Let $R$ be a Noetherian ring and $M$ a finitely generated $R$-module. Then, $M[[X]]$ is a torsion-free Auslander $R[[X]]$-module if and only if $M$ is a torsion-free Auslander $R$-module.

\begin{proof}
Since $M$ is finite and $R$ is Noetherian, $M$ is also a Noetherian $R$-module. This means that both the ring $R$ and the module $M$ have property (A). Now by Theorem \ref{FormalPowerAuslander} and Theorem \ref{FormalPowerZ}, the proof is complete.
\end{proof}
\end{theorem}

\section*{Acknowledgements}

The author was partly supported by the Department of Engineering Science at Golpayegan University of Technology and wishes to thank Professor Winfried Bruns for his invaluable advice. The author is also grateful for the useful comments by the anonymous referee.

\bibliographystyle{plain}

\end{document}